\def\C{{\bf C}}

\def\A{{\cal A}}
\def\H{{\cal H}}

\def\d{d}
\def\height{\mathop{\rm ht}}
\def\diam{\mathop{\rm diam}}

\def\span{\mathop{\rm span}}

\def\Co(#1,#2){C(#1,e^{-#2})}
\def\den#1{{\langle\gamma^{#1}1,1\rangle}}
\newcount\m
\def\t{\the\m\global\advance\m by 1}
\m=1
\newcount\n
\def\f{\the\n\global\advance\n by 1}
\n=1
\centerline{\bf Boundary unitary representations---}
\centerline{\bf right-angled hyperbolic buildings}
\medskip
\centerline{Uri Bader and Jan Dymara}
\footnote{}{UB: Partially supported by the ERC grant 306706. }
\footnote{}{JD: Partially supported by
Polish NCN grant
DEC-2012/06/A/ST1/00259.}

\bigskip
Abstract.
We study the unitary boundary representation
of a strongly transitive group acting on a right-angled
hyperbolic building.
We show its irreducibility.
We do so by associating to such a representation a representation of a certain Hecke algebra, which is a deformation of the classical representation of a hyperbolic reflection group.
We show that the associated Hecke algebra representation is irreducible.
\bigskip

\centerline{\bf 0. Introduction}
\medskip

Considering the group $SL(2,{\bf R})$ as the isometry group of the hyperbolic plane and the circle $S^1$ as its boundary,
one is led to study various unitary representations of $SL(2,{\bf R})$ on the function space $L^2(S^1)$: the so called ``principal series representations".
It is well known that these representations are irreducible, and constitute a sizable
part of the unitary dual of $SL(2,{\bf R})$.
This fundamental fact inspired many authors who obtained various generalizations.
For example, one may replace the field ${\bf R}$ with a non-archimedean local field, and the group $SL(2)$ with an arbitrary semi-simple algebraic group.
Important tools in the study of the representations thus obtained are the associated affine building and the so called Iwahori-Hecke algebra, see [Iw], [IM], [HM].
In another course of generalization, taken in [CS] (see also [FTP], [FTS], [Gar1], [BK]), the authors consider discrete subgroups
and show (for example) that the restriction of the principal $SL(2,{\bf R})$-representations on $L^2(S^1)$ to the subgroup $W$, generated by the reflections across the sides of a (compact) right angled polygon, is still irreducible.
Our current contribution links the previous two routes: we consider a group acting on a building and its unitary representation on $L^2$ of the boundary.
We use the aid of an associated Hecke algebra representation in order to analyze this boundary representation.
It turns out that the associated Hecke algebra representation could be seen as a deformation of the unitary representation of $W$ alluded to above.
We prove the irreducibility of this deformed representation using the tools developed in [BM] (see also the recent generalization [Gar2]).
Our exact theorem is:

\smallskip
\bf Theorem 0.1.\par\sl
Let $X$ be a right-angled (Lobachevsky) hyperbolic building
of finite thickness. Let $G$ be a group acting strongly transitively
on $X$. The associated unitary representation of $G$ on $L^2(\partial X)$
is irreducible.
\rm \par
\smallskip

As mentioned above, in the course of the proof of Theorem 0.1 we reduce it to the following theorem, which might be of independent interest.

\smallskip
\bf Theorem 0.2.\par\sl
Let $W$ be the group generated by the reflections across the codimension-1 faces of a compact right-angled polytope 
in the hyperbolic $n$-space.
Fix real parameters $q_s$, indexed by the codimension-1 faces of the polytope and consider the corresponding 
Iwahori-Hecke algebra $\H$.
Then the natural representation of $\H$ on $L^2(S^{n-1})$ is irreducible, provided that for every face $s$, $q_s\geq 1$.
\rm \par
\smallskip

The terms used in the formulations of the above theorems will be explained in the next sections.
In particular, the unitary representations considered are the principal series representation with the trivial parameter ($\epsilon=0$).
Unfortunately, the question whether other principal series representations are irreducible as well remains open.
The dimension of the buildings that we deal with is limited: Vinberg proved that compact 
right-angled hyperbolic polytopes do not exist in dimensions $>4$ (see [Dav, Cor.~6.11.7]).

The paper is divided into three parts. In part I we detail the setting of Theorem 0.1 and explain its reduction to Theorem 0.2.
Part II describes explicitly the Hecke algebra and its principal series representations.
Finally, in part III, we prove Theorem 0.2.

\bigskip
\centerline{\bf I. Representation on the boundary of the building}

\bigskip
\centerline{\bf 1. Definition of principal series}
\medskip

\def\a{1}
\m=1
\n=1

In this section we define a series of unitary representations of a
hyperbolic building automorphism group on the $L^2$-space of the
boundary of the building.

Let $P$ be a bounded polytope with finitely many faces in a hyperbolic space $\H^\d$.
Suppose that all dihedral angles of $P$ are of the form $\pi/k$, $k\in{\bf Z}$. Then
the set $S$ of reflections across the codimension-1 faces of $P$ generates a group $W$ acting on $\H^\d$,
called a \it hyperbolic Coxeter group. \rm
The action of $W$ on $\H^\d$ is geometric (cocompact and properly discontinuous).
The $W$-translates of $P$ have pairwise disjoint interiors and form a tessellation of $\H^\d$.

Let $X$ be a building whose Weyl group is a hyperbolic Coxeter group $W$, as above.
One can think of $X$ as of a set (of chambers) with a $W$-valued distance function.
One can also consider geometric realizations of $X$; one of them will be the hyperbolic realization $|X|$,
with chambers isometric to $P$ and apartments isometric to $\H^\d$ tessellated by copies of $P$ as above.
Codimension-1 faces of $P$ correspond to elements of $S$; this labeling extends consistently to
codimension-1 faces in $X$, the label of  such face usually called its \it type. \rm
In our tessellation of $\H^\d$, a codimension-1 face is shared by two chambers; in a building, there should be more of them, and their
number is called the \it thickness \rm of the building along the face. We assume that the thickness is finite and
depends only on the type $s$ of the face, and we denote it $q_s+1$. The \it thickness vector \rm $q=(q_s)_{s\in S}$ encodes
the thickness data. If all $q_s$ are equal, we interpret $q$ as their common value; this is the
\it uniform thickness \rm case. In this paper,
all formulae are written as if thickness was uniform. Yet, they can all
be easily re-interpreted to make sense for non-uniform thickness. Comments explaining some details
of this re-interpretation will be given in \S7.

An automorphism of $X$ is a bijection $X\to X$ preserving the $W$-valued distance. It can be realized geometrically
as an isometric map $|X|\to|X|$; that map preserves types of codimension-1 faces.
Finally, we want a group $G'<G=Aut(X)$ acting strongly transitively on $X$ (this means transitivity on the set of pairs (chamber $\in$ apartment)).
The existence of a (strongly) transitive action implies that the thickness of the
building along a face depends only on the type of that face.
The group $G$ can be equipped with the  compact--open topology coming from its action on $X$. We require $G'$ to be closed
in this topology (anyway, passing to closure preserves strong transitivity). Then both $G$ and $G'$ are locally compact, totally discontinuous, second countable,
generated by compact subgroups (eg by stabilizers of codimension-1 faces
of one chamber), unimodular. We normalize the Haar measure $\nu$  on $G$
by requiring that (every) chamber stabilizer has measure 1 (all such stabilizers are conjugate).

The metric space $|X|$ is $CAT(-1)$ and has a compact Gromov boundary $\partial |X|$, which we usually shorten to $\partial X$
(the set $X$ itself can be equipped with the gallery distance; then it is quasi-isometric to $|X|$, and has the same Gromov boundary).
The action of $G$ on $X$ extends to an action on $\partial X$. We will define a $G$-quasi-invariant
measure $\mu$ on $\partial X$ and consider the associated family of unitary $G$- and $G'$-representations
on $L^2(\partial X,\mu)$.

We fix an isometric identification of $\H^\d$ (along with $P$ and the tessellation
that it generates) with the
the Poincar\'e disc model ${\bf D}^\d$. We require that the center $0$ of the
model corresponds to an interior point of $P$, also to be denoted $0$.
The action map of $W$ to the orbit of $0$ is a quasi-isometry between
$W$ (with the $S$-word metric) and $\H^\d$ or ${\bf D}^\d$; we use it to identify
$\partial W$, $\partial \H^\d$ and $\partial{\bf D}^\d=S^{\d-1}$.
Hence, we get a measure $l$ on $\partial W$ corresponding to the
Lebesgue measure on the sphere.
For any chamber $c\in X$ and any apartment $A$ containing $c$ the $c$-based folding map identifies
$\partial A$ with $\partial W$; pulling $l$ back by this folding we get a measure $l_c$
on $\partial A$.
The stabilizer $G_c$ of $c$ in $G$ has a probability Haar measure $\nu_c$.
The map $p^c\colon G_c\times\partial A\ni(b,x)
\mapsto bx\in
\partial X$ is surjective by strong transitivity.
\smallskip
\bf Definition \a.\t\sl
$$\mu_c=p^c_*(\nu_c\times l_c);\quad \mu=\mu_{c_0}\,\, \hbox{\rm for some fixed chamber $c_0$}.\leqno(\a.\f)$$
\smallskip
\rm
Note that $\mu_c$ does not depend on the choice of $A$: any different $A'\ni c$ is of the form
$b'A$, for some $b'\in G_c$; then $(b,x)\mapsto (b(b')^{-1},b'x)$ transforms one variant of
$\nu_c\times l_c$ to the other while preserving fibres of $p^c$. On the other hand, $\mu_c$ does depend
on the choice of $c$.
\smallskip
\bf Lemma \a.\t \sl

The measures $\mu_c$, $\mu_{c'}$ are absolutely continuous with respect to each other for $c,c'\in X$.
\smallskip

\bf Lemma \a.\t \sl

For any $g\in G$ and any chamber $c\in X$ we have $g_*\mu_c=\mu_{g(c)}$.
\smallskip\rm

\rm In particular, for any $g\in G$ we get $g_*\mu=\mu_{gc_0}\!<\!\!<\!\mu$. Using this fact we now define
the principal series of $G$; Lemmas \a.2 and \a.3 will be proved later.
\smallskip
\bf Definition \a.\t\sl\par
\itemitem{a)} Let $X$ be a hyperbolic building of finite thickness,
and let $G$ be the group of all
type-preserving automorphisms of $X$. Assume that the action of $G$ on
$X$ is strongly transitive.
The principal series of $G$ is the family $\rho_\epsilon$ ($\epsilon\in
{\bf R}$) of representations of $G$ on $L^2(\partial X,\mu)$ given by
$$\rho_\epsilon(g)f\,(x)=f(g^{-1}x)
\left[ {d(g_*\mu)\over d\mu}(x)\right]^{{1\over2}+i\epsilon}.$$
\itemitem{b)} Let $G'$ be a closed subgroups of $G$, still acting strongly
transitively on $X$.
The principal series of $G'$ is the restriction of the above to $G'$.
\smallskip \rm
\bf Remark\rm\par
In the definition of $\rho_\epsilon$ one could replace
$\mu$ by any measure $\mu'$ satisfying $g_*\mu'\!<\!\!<\!\mu'$.
However, if $\mu'=h\mu$, $h\in L^1(\mu)$, then
the representation associated to $\mu'$ isometrically embeds in
that corresponding to $\mu$: the intertwiner is given by
$$L^2(\mu')\ni f\mapsto h^{{1\over2}+i\epsilon}f\in L^2(\mu).\leqno(\a.\f)$$
In particular, if we replace $c_0$
by a different chamber we get an equivalent representation.
%On the other hand, varying the shape of $P$ probably leads to non-equivalent representations.

\smallskip

\bf Proof of Lemma \a.3.
\par\rm Let  $A$ be an apartment containing $c$. We have a commutative diagram:

$$\matrix{
G_c\times\partial A&\mathop{\longrightarrow}\limits^{\iota_g\times g}&G_{g(c)}\times\partial(gA)\cr
\Big\downarrow\rlap{$\vcenter{\hbox{$\scriptstyle p^c$}}$}
& &\Big\downarrow\rlap{$\vcenter{\hbox{$\scriptstyle p^{g(c)}$}}$}  \cr
\partial X&\mathop{\longrightarrow}\limits^{g}&\partial X}\leqno(\a.\f)$$
where $\iota_g$ denotes conjugation by $g$. Now:
$$g_*(\mu_c)=g_*p^c_*(\nu_c\times l_c)=p^{g(c)}_*(\iota_g\times g)_*(\nu_c\times l_c)=
p^{g(c)}_*(\nu_{g(c)}\times l_{g(c)})=\mu_{g(c)}.\leqno(\a.\f)$$
\hfill{$\diamond$(Lemma \a.3)}
\smallskip
\bf Proof of Lemma \a.2. \par\rm
We may assume that $c$ and $c'$ are $s$-adjacent for some given $s$ (we get the
general statement by applying this special case along a gallery). 
For convenience, we will assume $c=c_0$ (as its choice was arbitrary) and $c'=c_1$,
where 
$c_1,\ldots c_q$ denote the
$s$-neighbors of $c_0$. Let $X_i$ be the set of chambers that are closer to $c_i$
than to any other $c_j$. By $\partial X_i$ we denote the set of points in $\partial X$ that are in the closure
of $X_i$. We put $\mu_i=\mu_{c_i}$, $p^i=p^{c_i}$, $l_i=l_{c_i}$, $G_i=G_{c_i}$. We denote by
$\pi_j$ the $c_j$-based folding map $X\to W$, as well as its extension $\partial X\to\partial W$.
For $x\in\partial W$ put $r(x)={d(s_*l)\over dl}(x)$.

\smallskip

\bf Lemma \a.\t\sl\par
Fix an $i>0$, and let $q=q_s$. Then
$$\mu_i|_{\partial X_i}=q(r\circ\pi_0)\mu_0|_{\partial X_i}.\leqno(\a.\f)$$
\rm\par Proof.
We choose some apartment $A$ containing $c_0$ and $c_i$. Then $\partial A$ will be the support of
both $l_0$ and $l_i$.
We put $\partial A_0=\partial A\cap\partial X_0$, $\partial A_i=\partial A\cap\partial X_i$.
We will prove the following
sequence of equalities:
$$\eqalign{\mu_i|_{\partial X_i}&=p^i_*(\nu_{G_i}\times l_i)|_{\partial X_i}=
p^i_*(q\nu_{G_i\cap G_0}\times l_i)|_{\partial X_i}\cr &=
p^0_*(q\nu_{G_i\cap G_0}\times l_i)|_{\partial X_i}=
q(r\circ\pi_0)p^0_*(\nu_{G_0}\times l_0)|_{\partial X_i}=
q(r\circ\pi_0)\mu_0|_{\partial X_i}.}\leqno(\a.\f)$$
\par
$\bullet$ The first and last equality are immediate consequences of the
definitions of $\mu_i$ and $\mu_0$.\par
$\bullet$ $2^{\rm nd}$ equality: We have
$(\bigcup_{j\ne 0}\partial X_j)/G_0=\partial A_i$. Moreover, the setwise stabilizer of
$\partial X_i$ in $G_0$ is $G_0\cap G_i$.
Therefore $\partial X_i/(G_0\cap G_i)=\partial A_i=
\partial X_i/G_i$.
It follows that for cosets $(G_i\cap G_0)b$, $b\in G_i$, we  have
$(G_i\cap G_0)b\,\partial A_i=\partial X_i$. In particular, for each
$x\in\partial A_i$ we can choose a system $b_j$, $j=1,2,\ldots,q$ of
representatives of right cosets of $G_i\cap G_0$ in $G_i$ such that
$b_jx=x$.

Let us pick a measurable $U\subseteq\partial X_i$. Fix an $x\in\partial A_i$
and a system $b_j$ as above (for this $x$).
Then
$$\{b\in G_i\mid bx\in U\}=\bigcup_j\{b\in G_i\cap G_0\mid bx\in U\}b_j.\leqno(\a.\f)$$
Since right multiplication by $b_j$ preserves the Haar measure,
we have
$$\nu(\{b\in G_i\mid bx\in U\})=q\nu(\{b\in G_i\cap G_0\mid bx\in U\}).\leqno(\a.\f)$$
Integrating this equality over $\partial A_i$ with respect to $l_i$
we get that $p^i_*(\nu|_{G_i}\times l_i)$ and $p^i_*(q\nu|_{G_i\cap G_0}\times
l_i)$ agree on $U$.
\par
$\bullet$ $3^{\rm rd}$ equality: The pushed measure is supported on the
intersection of the domains of $p^0$ and $p^i$; both maps are defined
by the same formula.
\par
$\bullet$ $4^{\rm th}$ equality: we have $s_*l=rl$, and (on $\partial A$) $\pi_i=s\circ\pi_0$.
In the following calculation the domains of $\pi_0$ and $\pi_i$ are restricted to $\partial A$:
$$l_i=(\pi_i^{-1})_*l=((s\circ\pi_0)^{-1})_*l=(\pi_0^{-1})_*s_*l=(\pi_0^{-1})_*(rl)=
(r\circ\pi_0)(\pi_0^{-1})_*l=(r\circ\pi_0)l_0.\leqno(\a.\f)$$
On $G_i\cap G_0$ the measures $\nu|_{G_0}$ and
$\nu|_{G_i\cap G_0}$ coincide;  $p^0$ maps the remaining part of
$G_0$ ($\times\partial A_i$) outside $\partial X_i$.

\hfill{$\diamond$(Lemma \a.5)}
\smallskip
To finish the proof of Lemma \a.2 we use principles of symmetry. First
$$\mu_i|_{\partial X_0}={1\over q(r\circ\pi_i)}q(r\circ\pi_i)\mu_i|_{\partial X_0}
={1\over q(r\circ\pi_i)}\mu_0|_{\partial X_0}
={1\over q(r\circ s\circ\pi_0)}\mu_0|_{\partial X_0}
=q^{-1}(r\circ \pi_0)\mu_0|_{\partial X_0},
\leqno(\a.\f)$$
the last equalities due to: $\pi_i=s\circ\pi_0$ on $\partial X_0$; $r\circ s={1\over r}$.
Next, for $j\ne0,i$, we apply (\a.10) twice:
$$\mu_i|_{\partial X_j}=q^{-1}(r\circ \pi_j)\mu_j|_{\partial X_j}=\mu_0|_{\partial X_j}.\leqno(\a.\f)$$
\hfill{$\diamond$(Lemma \a.2)}
\smallskip
We combine (\a.5), (\a.10) and (\a.11) into the following corollary
(recall that $\mu=\mu_0$ and $\pi=\pi_0$).
\smallskip
\bf Corollary \a.\t\sl

Suppose $\phi\in G$ and $\phi c_0=c_i$ ($i>0$). Then
$$\phi_*\mu=\cases{
\mu& on $\partial X_j$, $j\ne0,i$;\cr
q(r\circ\pi)\mu& on $\partial X_i$;\cr
q^{-1}(r\circ\pi)\mu& on $\partial X_0$.
}\leqno(\a.\f)$$
\smallskip\rm
\bf Remark \a.\t.\rm\par
It is possible to deduce from Corollary \a.6 a general formula for ${d\mu_c\over d\mu_{c'}}$
in terms of Busemann functions. In a metric space $(X,\rho)$ one defines the Gromov
product $(x|y)_z={1\over2}(\rho(z,x)+\rho(z,y)-\rho(x,y))$. If the space is hyperbolic
and has boundary $\partial X$, we may try to extend the Gromov product by putting
$(x|b)_z=\lim_{y\to b}(x|y)_z$. The limit is taken over a sequence of $y\in X$ converging
to $b\in\partial X$; the limit may not exist, or may depend on the choice of the sequence of $y$'s,
but in our cases of interest these problems will happen only for $b$ in some zero-measure set. Finally, the Busemann function is
$\beta^\rho_b(x,y)=\rho(x,y)-2(y|b)_x$ (for $x,y\in X$, $b\in\partial X$).
In our case, there are two metrics of interest on the building $X$.

1) Each chamber $c\in X$ has a geometric realization $|c|$ in $|X|$, canonically
isometric to the fundamental polytope $P$. Inside $|c|$ there is a copy $0_c$ of the
point $0\in P$. We put $|c-c'|=|0_c-0_{c'}|$ (the distance from $0_c$ to $0_{c'}$ in $|X|$).

2) The gallery distance: $\ell(c,c')$ is the length of a shortest gallery in $X$
starting at $c$ and ending at $c'$.

Now we can state the formula:
$${d\mu_c\over d\mu_{c'}}=e^{-(\d-1)\beta^{|\cdot|}_b(c',c)}q^{-\beta^\ell_b(c',c)}.\leqno(\a.\f)$$
For $c$ adjacent to $c'$ this formula reduces to Corollary \a.6. The general case follows
from this special case and the cocycle property of a Busemann function:
$$\beta^\rho_b(x,z)=\beta^\rho_b(x,y)+\beta^\rho_b(y,z).\leqno(\a.\f)$$

As noticed by Garncarek, the right hand side of (\a.13) can be rewritten as
$$\exp(-(\d-1)\beta^{|\cdot|}_b(c',c)-\ln(q)\beta^\ell_b(c',c))=
\exp(-(\d-1)\beta^{\rm mix}_b(c',c)),\leqno(\a.\f)$$ where
${\rm mix}(c',c)=|c'-c|+{\ln{q}\over \d-1}\ell(c',c)$.
This means that the measures $(\mu_c)_{c\in X}$ form a quasi-conformal
family with respect to the mix metric. Thus, the results of [Gar2] apply
to our setting, yielding (in many cases) a different proof of Theorem 0.1.
Garncarek requires the group to be discrete; his result can be applied
to a cocompact lattice in $G'$. In the right-angled case, such a lattice
exists in the
automorphism group $G$ and in many interesting subgroups $G'$, cf [CT].

This remark holds also in the multi-parameter case. One replaces the gallery
metric $\ell$ by a family $(\ell_s)_{s\in S}$ of pseudo-metrics:
$\ell_s(c,d)$ is the number of $s$-type codimension-1 faces traversed by
a minimal gallery from $c$ to $c'$. Then the mix metric  is defined as
${\rm mix}(c',c)=|c'-c|+\sum_{s\in S}{\ln{q_s}\over \d-1}\ell_s(c',c)$.

\bigskip

\centerline{\bf 2. General nonsense}
\medskip

\def\a{2}
\m=1
\n=1
Suppose that $(V,\rho)$ is a unitary representation of a locally compact group $G$,
and let $K$ be a compact-open subgroup of $G$. Then one defines the Hecke
algebra $H(G,K)$ as the convolution algebra of all compactly supported $K$-bi-invariant
functions on $G$. Elements of $H(G,K)$ are continuous (even locally constant),
because $K$ is open. They act on $V$ by
$$\rho(f)v=\int_Gf(g)\rho(g)v\,dg.\leqno(\a.\f)$$
This action preserves the space of $K$-invariants $V^K=\{v\in V\mid (\forall k\in K)(\rho(k)v=v)\}$,
due to left $K$-invariance of $f$.

\smallskip
\bf Proposition \a.\t.\sl\par
Suppose that:
\item{1)}{ $V^K$ is $G$-cyclic in $V$ (ie $V$ is the closure of the linear span of $\rho(G) V^K$);}
\item{2)}{ $V^K$ is Hecke-irreducible (ie $V^K$ has no non-trivial closed $H(G,K)$-invariant subspace).}

Then $V$ is an irreducible $G$-representation.
\rm\par
Proof. Suppose not. Let $V=V_0\oplus V_1$ be a non-trivial orthogonal decomposition
into subrepresentations. Then $V^K=V_0^K\oplus V_1^K$. The subspaces $V_0^K$ and $V_1^K$ are
non-zero: if $V_0^K$ was zero, we would have $\rho(G)V^K=\rho(G)V_1^K\subseteq V_1$,
contradicting the cyclicity assumption; similarly for $V_1^K$. But this contradicts
Hecke-irreducibility of $V^K$.
\hfill{$\diamond$}

The goal of Part I is to establish $G$-cyclicity of $V^K$ in our context. To do that,
we need to recall certain facts about hyperbolic right-angled buildings and their
automorphism groups.
\bigskip

\centerline{\bf 3. Right-angled buildings}
\medskip

\def\a{3}
\m=1
\n=1

A Coxeter group $(W,S)$  is \it right-angled \rm if any two elements of $S$
either commute or span an infinite dihedral subgroup of $W$. A building is
right-angled if its Weyl group is right-angled. A hyperbolic reflection group
associated to a polyhedron $P$ is right-angled if all dihedral angles of $P$
are $\pi/2$.

Morphisms of right-angled buildings are discussed at length in section 4 of [DO].
Here we just summaries the results.
Let $X$ be a right-angled building. We fix a chamber $c$.
A set $Y\subseteq X$ is \it star-like \rm
(with respect to $c$), if every minimal gallery from $c$ to any $y\in Y$ is contained in $Y$.
We choose a well-ordering $<$ on $X$ such that all initial segments $X_{<x}$ are star-like.
A morphism $\phi\colon X\to X$ can be constructed inductively with respect to $<$. We first arbitrarily
choose $\phi(c)$. If $\phi$ is defined on $X_{<x}$, we try extend it  to $x$. Two things can happen:
\item{1)} (freedom) All minimal galleries from $x$ to $c$ are of the form $(x,y,\ldots,c)$
for some unique $y$. Then there is choice: if $s$ is the type of the common face of $x$ and $y$, then
$\phi(x)$ can be chosen to be any of the $s$-neighbors of $\phi(y)$.
\item{2)} (determinism) There are at least two minimal galleries from $x$ to $c$ starting differently:
$(x,y,\ldots,c)$ and $(x,y',\ldots,c)$ with $y\ne y'$. Then there is a unique choice of $\phi(x)$
that extends (the so-far constructed part of) $\phi$ to a morphism.
The values of $\phi(y)$, $\phi(y')$, and similar values given by
other galleries consistently and uniquely determine $\phi(x)$.

\noindent Whether $x$ falls into 1) or 2) depends only on $w=\pi_c(x)\in W$. We compare
the $S$-word lengths $\ell(w)$ and $\ell(ws)$: if $\ell(w)>\ell(ws)$ for just one $s\in S$,
then we are in case 1); otherwise we are in case 2). In [DO] the set of all pairs $(y,s)$ as in
case 1) is called the \it root set \rm of $X$ and denoted $R(X)$. To check whether the constructed morphism
is an automorphism we have the following criterion: for each $(y,s)\in R(X)$ the map
$\phi$ restricts to a bijection between the sets of $s$-neighbors of $y$ and of $\phi(y)$.
In particular, if the thickness of $X$ is type-dependent (given by a thickness vector),
then a partial automorphism defined on an initial segment of $X$ can always be extended to an automorphism.

To make use of the above procedure it is necessary to have well-orderings with star-like
initial segments (let us temporarily call them \it nice\rm). For buildings of finite thickness
nice orderings are not hard to come by, here are some examples.
\item{-} Any ordering compatible with gallery distance from $c$ (ie satisfying
$\delta(c,x)<\delta(c,y)\Rightarrow x<y$) is nice.
\item{-} For star-like subsets $A,B$ of $X$, we can put a nice ordering on $A$,
then extend it to $A\cup B$ (so that $A$ is an initial segment), and then extend it to $X$
(keeping $A\cup B$ as an initial segment). We will use this type of ordering for
convex sets $A$, $B$ containing $c$.
\item{-} Any refinement of a nice ordering of $W$ is nice. By a refinement we mean
an ordering on $X$ satisfying $x<y\iff \pi_c(x)<\pi_c(y)$.
\smallskip
Finally, we need to discuss \it standard open neighborhoods \rm (cf [DO], section 2).
Let $X$ be a hyperbolic right-angled building with finite thickness and Weyl group $W$.
We distinguish a chamber $c_0$ and denote by $\pi$ the $c_0$-based folding map.
Pick any wall $H$ in $|W|$. Let $\overline{H}$ is the completed wall, ie the closure of
$H$ in $|W|\cup\partial W$. Consider the connected components
of $(|X|\cup\partial X)\setminus \pi^{-1}(\overline{H})$. The ones that do not contain $c_0$ are called
\it standard (open) neighborhoods. \rm Every point $p\in \partial X$ has a basis of open
neighborhoods in $|X|\cup\partial X$ consisting of suitably chosen standard neighborhoods.
The intersections of sets from this basis with $\partial X$ form a basis of neighborhoods of
$p$ in $\partial X$.

Any standard open neighborhood $U$ contains a unique minimal chamber: a chamber  $x$ in $U$  with minimal
gallery distance to $c_0$. This $x$ satisfies the freedom condition 1).  Conversely, any chamber $x$ satisfying 1)
is the minimal chamber of a unique standard neighborhood, to be denoted $U(x)\cup\partial U(x)$
($U(x)$ for the part in $|X|$, $\partial U(x)$ for the part in $\partial X$).

\bigskip
\centerline{\bf 4. Invariants are cyclic}
\def\a{4}
\m=1
\n=1
\medskip
The standing notation throughout the paper is as follows:
$X$ is a right-angled hyperbolic building of finite thickness
$q$ (either uniform or multi-parameter), with Gromov boundary $\partial X$;
$G$ is the type-preserving automorphism group of $X$, equipped with the compact-open topology;
$c_0$ is a fixed chamber in $X$ (called ``the base chamber''),
and $K$ is the stabilizer of $c_0$ in $G$;
$G'$ is a closed subgroup of $G$, acting strongly transitively on $X$;
we put $K'=G'\cap K$;
$(V,\rho_\epsilon)$ is a principal series representation of $G$,
or its restriction to $G'$.

In this section we prove that $V^{K'}$ is $G'$-cyclic in $V$.
We first reduce to the case of $G$ and $K$.
\smallskip

\bf Lemma \a.\t. \sl \par
\item{a)} $V^K=V^{K'}$. More generally, $V^{G_c}=V^{G'_c}$ for any chamber $c\in X$.
\item{b)} $\span\rho_\epsilon(G)V^K=\sum_{c\in X}V^{G_c}=\span\rho_\epsilon(G')V^{K'}$.
\rm\par
Proof.
We have $\partial X/K\simeq \partial W\simeq \partial X/K'$, the quotient maps being $\pi$ in both cases.
Since $\pi$ is measure preserving, we get $V^K\simeq L^2(\partial W,l)\simeq V^{K'}$ as Hilbert spaces,
the isomorphisms given by
precomposition with $\pi$.

For $g\in G$ we have $\rho_\epsilon(g)V^K=V^{gKg^{-1}}=V^{G_{gc_0}}$. To finish the proof of a),
let us choose a $g\in G'$, such
that $gc_0=c$; we then have $V^{G_c}=\rho_\epsilon(g)V^K=\rho_\epsilon(g)V^{K'}=V^{G'_c}$.
Part b) follows as well:
$$\span\rho_\epsilon(G)V^K=\sum_{g\in G}V^{G_{gc_0}}=\sum_{c\in X}V^{G_c}=\sum_{c\in X}V^{G'_c}=
\sum_{g\in G'}V^{G'_{gc_0}}=\span\rho_\epsilon(G')V^{K'}\leqno(\a.\f)$$
\hfill{$\diamond$}
\smallskip
\bf Theorem \a.\t.\par\sl
The linear span of $\rho_\epsilon(G)V^K$ is dense in $V$.\rm\par
Proof.
Assume, by contradiction, that $\span\rho_\epsilon(G)V^K=\sum_{c\in X}V^{G_c}$ is not dense in $V$.
Then there
exists a non-zero function $f\in \left(\sum_{c\in X}V^{G_c}\right)^\perp$. The latter space being a (closed)
$G$-subrepresentation, we may average $f$ over a small subgroup of $G$ and still get a non-zero $f$ in the same
space; in other words, we may assume that $f$ is invariant under $G_Y$---the pointwise stabilizer of a sufficiently
large finite subset $Y\subseteq X$. We may assume and do assume that $c_0\in Y$.

It turns out that, for diligently chosen chamber $c\in X$ and open set $U\subseteq \partial X$,
the properties of $G_Y$-invariance and $G_c$-invariance are equivalent on $U$.
Then, the function $f\chi_U$, if non-zero, will be $G_c$-invariant and not perpendicular to $f$,
yielding a contradiction ($\chi_U$ is the characteristic function of $U$).
We proceed to the details.
\smallskip
Let the gallery diameter of $Y$ be $M$.
Consider the gallery-distance $M$-ball around $1$ in $W$, and
the finitely many  tessellation walls $H'$ in $|W|$ that intersect (the closed geometric realization of) this ball.
The union of those walls:
${\cal W}=\bigcup\{\overline{H}\mid \hbox{\rm gallery distance
from $1$ to $H'$ is} \le M\}$, is closed in $|W|\cup\partial W$.
Also, ${\cal W}\cap \partial\H^\d$ has zero Lebesgue measure---this set is a finite union
of zero-measure boundaries of walls. We deduce that $\pi^{-1}(\partial{\cal W})$ is closed in $\partial X$
and has $\mu$-measure zero. It follows that there exists $x\in\partial X\setminus\pi^{-1}(\partial{\cal W})$
such that $f$ has non-zero restriction to every open neighborhood of $x$. Let us pick a standard open neighborhood
$U(c)\cup\partial U(c)$ of $x$ in $|X|\cup\partial X$ disjoint from  the closed set $\pi^{-1}({\cal W})$.
\smallskip
\bf Lemma \a.\t. \sl \par
For every $g\in G_c$ there exists $g'\in G_c\cap G_Y$ such that $g|_{U(c)}=g'|_{U(c)}$.
\rm\par

This lemma implies Theorem \a.2, as follows.
The set $\partial U(c)$ is $G_c$-invariant; every $g\in G_c$ preserves $\mu|_{\partial U(c)}$ and
$f\chi_{\partial U(c)}$, since so does $g'$ given by the lemma.
It follows that
$f\chi_{\partial U(c)}\in V^{G_c}$, $\langle f,f\chi_{\partial U(c)}\rangle=\int_{\partial U(c)}|f|^2\,d\mu>0$---contradiction.

Proof of Lemma \a.4: We choose an apartment $A$ containing $c_0$ and $c$. The map $\pi\colon A\to W$
is an isomorphism; we can therefore move all the data (the gallery $M$-ball around 1, the set ${\cal W}$)
from $W$ to $A$.
In this proof it will be convenient to consider foldings onto $A$ (rather than $W$). Thus, $\pi$ will denote
the folding map $X\to A$ fixing $c_0$, and $\pi_c$ will fix $c$.
Let $H$ be the boundary wall of $\pi(U(c))$. Cut $|A|$ along the walls $H'$ that intersect $H$;
let $C$ be the component of $c_0$ in $|A|\setminus\bigcup\{H'\mid H'\cap H\ne \emptyset\}$.
By the definition of ${\cal W}$, the gallery $M$-ball around $c_0$ in $A$ is contained in $C$; since
$\pi_c$ is a retraction, $Y\subseteq \pi_c^{-1}(C)$. We choose a well-ordering on $X$ starting in $c$, so that:
chambers in $U(c)$ form an initial segment; chambers in $\pi_c^{-1}(C)$ form the next segment; all initial segments are star-like.
Then start defining the automorphism $g'$ by imposing: $g'=g$ on $U(c)$;  $g'={\rm id}$ on $\pi_c^{-1}(C)$.
This defines a partial automorphism on an initial segment, which extends to an automorphism of $X$.
\hfill{$\diamond(\hbox{\rm Lemma \a.4, Theorem \a.2})$}
\bigskip

\centerline{\bf II. Right-angled Hecke algebra}
\bigskip
Our goal in this part is to define and discuss the principal series representations of the Hecke algebra on $V^K\simeq L^2(\partial W)$.
Elements of the Hecke algebra $H(G,K)$ are compactly supported functions on $G$ that are
$K$--bi-invariant. In other words, they correspond to finitely supported functions
on $K\backslash G/K$, which, by the Bruhat decomposition, is naturally identified with $W$.
Thus, as vector spaces, $H(G,K)\simeq \C[W]$. The convolution multiplication in $H(G,K)$
does not correspond to the group algebra multiplication in $\C[W]$, but rather
to its deformation described below in (5.1). The isomorphism is explained in
[Bourb, Ch.~IV, \S 2, ex.~24] and in  [Dav, Lemma 19.1.5]. The source
[Dav] proves the result not only for the full automorphism group $G$, but also
for its strongly transitive subgroups $G'$.
In \S5 below we describe the Hecke algebra structure on $\C[W]$ associated with an abstract choice of parameter $q$.
In \S6 we study the representation of $H(G,K)$ on $V^K$.
In both sections we consider, only for the simplicity of the discussion, the uniform thickness case $q$.
In \S7 we will explain how to read through the previous sections while considering the non-uniform thickness case.
We will also explain how to define the principal series representations for an arbitrary Hecke algebra, that is when considering complex parameters $q_s$.

\bigskip
\centerline{\bf 5. Multiplication in right-angled Hecke algebras}
\def\a{5}
\m=1
\n=1
\medskip
In this section $(W,S)$ is a right-angled Coxeter group, not necessarily hyperbolic.
The Hecke algebra of $W$ is the space
$\H=\oplus_{w\in W}\C e_w$  with $\C$-bilinear multiplication
given by
$$e_se_w=\cases
{e_{sw} &if $\ell(sw)>\ell(w)$\cr
(q-1)e_w+qe_{sw}&if $\ell(sw)<\ell(w)$,
}\leqno(\a.1)
$$
where $w\in W$ and $s\in S$.
We can rewrite this as a formula for $e_sf\,(u)$, where
$f\in \H=\C[W]$ is treated as a finitely--supported function on $W$:
$$ e_sf\,(u)=\cases
{qf(su)&if $\ell(su)>\ell(u)$\cr
(q-1)f(u)+f(su)&if $\ell(su)<\ell(u)$.
}\leqno(\a.2)
$$
Our first goal is to establish a similar formula for $e_wf$ without assuming
that $w$ is a generator. The formula will hold for every right-angled Coxeter group,
regardless of its hyperbolicity.

\bf Definition \a.\t \rm
\item{a)}
For $A,D\subseteq W$ we denote by $P(A|D)$ the set of walls in $W$
separating $A$ from $D$. In $|W|$ a wall divides $|W|$ into two connected components;
this division induces a partition of the set of chambers, hence of $W$, into
two parts. This partition of $W$ is the combinatorial meaning of ``a wall in $W$''.
We shall often identify walls with the corresponding reflections.
In particular, for $s\in S$ and  a wall/reflection $a$ we denote by $a^s$ the wall $s(a)$/the reflection $sas$;
for sets we put $A^s=\{a^s\mid a\in A\}$.
There is a natural poset structure on $P(A|D)$: $H\le H'$ if the half space bounded by $H$ containing $A$
is a subset of the half space bounded by $H'$ containing $A$.
\item{b)}
By $\A P(A|D)$ we denote the set of anti-chains in
$P(A|D)$. Notice that if $h\in\A P(A|D)$, than any two elements of $h$ commute; the converse is also true: a pairwise commuting family of wall reflections in $P(A|D)$ is an
anti-chain in $P(A|D)$.
We denote by $h$ the product of all elements of $h$
(the context will make the meaning clear). Observe that $h^{-1}=h$.
\item{c)}
For $h\in \A P(A|D)$ we define an ``interval'' $[A,h]$ as the set of all
$H'\in P(A|D)$ such that no $H\in h$ is smaller than $H'$:
$[A,h]=\{H'\in P(A|D)\mid (\forall H\in h)\,\neg (H<H')\}$ (notice that
$[A,\emptyset]=P(A|D)$). We define the height of $h$ as $\height(h)=\# [A,h]$.
\medskip

\bf Proposition \a.1\sl \par In the Hecke algebra of a right-angled Coxeter group
$$e_wf\,(u)=\sum_{h\in\A P(1|u,w)}(q-1)^{\# h}q^{\ell(w)-\height(h)}f(w^{-1}hu).\leqno(\a.3)$$
\par\rm
Proof. Induction on $\ell(w)$. The length $1$ case is covered by (\a.2).
Now assume that the formula is true for $w$ and that $\ell(sw)>\ell(w)$; we will deduce that the formula holds for $sw$.
There are two cases: either $\ell(su)>\ell(u)$, or $\ell(su)<\ell(u)$.

\smallskip
Case 1. $\ell(su)>\ell(u)$.

Then $1$, $w$ and $u$ are on the same side of the wall $s$.
We have
$$
\eqalign{
e_{sw}f\,(u)&=e_s(e_wf)\,(u)=q(e_wf)\,(su)\cr
&=q\sum_{h\in\A P(1|su,w)}(q-1)^{\# h}q^{\ell(w)-\height(h)}f(w^{-1}hsu)
}\leqno(\a.4)
$$
while the postulated formula for $e_{sw}f\,(u)$ is
$$\sum_{h\in\A P(1|u,sw)}(q-1)^{\# h}q^{\ell(sw)-\height(h)}f(w^{-1}shu).\leqno(\a.5)$$
The equality between (\a.4) and (\a.5) follows from three observations:
(1) $\ell(sw)=\ell(w)+1$; (2) $P(1|u,sw)=P(1|su,w)$; (3) if $c\in P(1|su,w)$ then $sc=cs$.
The first of them is clear. To prove the other two we need a lemma.
\smallskip
\bf Lemma \a.2\par\sl
Let $s\in S$; suppose $s\in P(1,g|h)$ and $c\in P(1|h,g)$. Then $sc=cs$.\par
\rm Proof.
For a wall $H$ we denote by
$H^-$  the half--space bounded by $H$ that contains $1$; $H^+$ is the other half--space
bounded by $H$.
If $sc\ne cs$, the walls $s$ and $c$ do not intersect. Consequently, one of the following holds:
$s^+\subseteq c^+$, $s^-\subseteq c^+$, $s^+\subseteq c^-$, $s^-\subseteq c^-$.
However,
$s\in s^+\cap c^-$,
$1\in s^-\cap c^-$,
$h\in s^+\cap c^+$,
$g\in s^-\cap c^+$.
Hence, neither inclusion may hold.
\hfill{$\diamond$}

Putting $g=w$ and $h=su$ we deduce observation (3) from Lemma \a.2.
Finally, $P(1|su,w)=P(1|su,w)^s=P(s|u,sw)=P(1|u,sw)$:
the first equality follows from (3); the second is clear if we interpret elements of $P(1|su,w)$ as walls;
the third is true because $s$ (the only wall separating elements $1$ and $s$) does
not belong to either side.
\smallskip
Case 2. $\ell(su)<\ell(u)$.

We have
$$\eqalign{
e_{sw}f\,(u)=&e_s(e_wf)\,(u)=(q-1)(e_wf)\,(u)+e_wf\,(su)\cr
=&(q-1)\sum_{h\in\A P(1|u,w)}(q-1)^{\# h}q^{\ell(w)-\height(h)}f(w^{-1}hu)\cr
&+
\sum_{h\in\A P(1|su,w)}(q-1)^{\# h}q^{\ell(w)-\height(h)}f(w^{-1}hsu),
}\leqno(\a.6)
$$
which we need to compare with
$$
\sum_{k\in\A P(1|u,sw)}(q-1)^{\# k}q^{\ell(sw)-\height(k)}f(w^{-1}sku).\leqno(\a.7)$$
We split the sum in (\a.7) into two parts: over $k\ni s$ and over $k\not\ni s$.
We will show that these  parts are equal to the two sums in (\a.6).
\smallskip
\bf Lemma \a.3\par\sl
In $P(1|u,sw)$ we have $[1,s]=P(1|u,w)\cup\{s\}$.
\par\rm
Proof. Suppose $c\in[1,s]$ and $c\ne s$. Then $c\in P(1|u,sw)$ implies $c\in P(1|su,w)$;
therefore $c\in P(1|u,w)$.
Suppose now that $c\in P(1|u,w)$. Lemma \a.2 (with $g=w$, $h=u$) implies that
$cs=sc$, hence $c\in[1,s]$. \hfill{$\diamond$}

It follows from this lemma that $h\mapsto k=h\cup\{s\}$ is a bijection between
$\A P(1|u,w)$ and $\{k\in \A P(1|u,sw)\mid k\ni s\}$. We have
$\# k=\# h+1$, $\height(k)=\height(h)+1$, $sk=h$, $\ell(sw)=\ell(w)+1$. Thus, the first sum in (\a.6) is equal to
the first sub-sum of (\a.7).

Notice that $P(1|u,sw)^s=P(s|su,w)=P(1|su,w)\cup\{s\}$ (disjoint union).
If $k\not\ni s$, then we put $h=k^s\in P(1|su,w)$. We have $\# k=\#h$, $sk=hs$,
$\height(k)=\height(h)+1$ (for $[1,h]=[1,k]^s-\{s\}$). Therefore, the second sum in (\a.6)
is equal to the second sub-sum of (\a.7).
\hfill{$\diamond$}

\bigskip
\centerline{\bf 6. Hecke algebra action in principal series representations}
\medskip
\def\a{6}
\m=1
\n=1

In this section we will describe the Hecke algebra representation on $V^{K'}$ for $V$ in the
principal series of $G'$
(the description is the same for $G$ and for $G'$).
In  this case $V=L^2(\partial X,\mu)$ and $V^{K'}=L^2(\partial X,\mu)^{K'}\simeq
L^2(\partial W,l)$, the last isomorphism
given by composition with $\pi$. Thus all Hecke representations (for all $q$) are
realized on the same Hilbert space
$L^2(\partial W,l)$. This includes the case $q=1$, when we get the principal series representation of the group
algebra of $W$. The dependence on $q$ being important for us, we should denote the action of $e_w$ on $f$ by
$\rho_\epsilon^q(e_w)f$; to shorten this we will write $w^qf$ instead.

Our goal for the remainder of this section
is to express the Hecke actions $w^qf$ in terms of the Weyl actions $w^1f$. We start with the case $w=s\in S$.
The wall $s$ separates $\partial W$ into two pieces: $\partial s^-$
(on the side of $1$) and
$\partial s^+$ (on the side of $s$).
The next lemma is somewhat analogous to (5.2).

\smallskip
\bf Lemma \a.\t\sl\par
Suppose $f\in L^2(\partial W,l)$. Then, for $x\in\partial W$:
$$(s^qf)(x)=\cases{
q^{{1\over2}-i\epsilon}(s^1f)(x)& for $x\in \partial s^-$\cr
(q-1)f(x)+ q^{{1\over2}+i\epsilon}(s^1f)(x)& for $x\in \partial s^+$
}\leqno(\a.\f)$$
\par
\rm Proof.
We pick some apartment $A$ containing the base chamber $c_0$
and its $s$-neighbor $c_1$.
Let $c_2\ldots,c_q$ be the other $s$-neighbors of $c_0$.
Put $F=f\circ\pi$. Let $y$ be the unique point in $\partial A$ that is mapped
to $x$ by $\pi$.
Let $K'_s$ be the stabilizer (in $G'$) of the common face of $c_0$ and $c_1$.
Under the isomorphism $\H\simeq\C[K'\backslash G'/K']$ the element $e_s$ corresponds to
the indicator function of (the set-theoretic difference) $K'_s\setminus K'$.
For $i>0$ choose $s_i\in G'$ that exchanges $c_0$ and $c_i$.
Then we have a disjoint union decomposition $K'_s\setminus K'=\bigcup_i s_iK'$.
$$\eqalign{
(s^qf)(x)
=\int_{K'_s\setminus K'}\rho_\epsilon(g)F\,(y)\,dg
&=\sum_i\int_{K'}F(b^{-1}s_i^{-1}y){\left[{d(s_ib)_*\mu\over d\mu}(y)\right]}^{{1\over2}+i\epsilon}db\cr
&=\sum_i\int_{K'}F(s_i^{-1}y){\left[{d(s_i)_*\mu\over d\mu}(y)\right]}^{{1\over2}+i\epsilon}db}
\leqno(\a.\f)$$
The measure $(s_i)_*\mu$ has been calculated in Cor.~1.6.
Suppose $x\in \partial s^-$. Then $F(s_i^{-1}y)=f(sx)$ and $((s_i)_*\mu)(y)=\mu_i(y)=q^{-1}r(x)\mu(y)$.
Therefore
$$(s^qf)(x)
=qf(sx)q^{-{1\over2}-i\epsilon}r(x)^{{1\over2}+i\epsilon}
=q^{{1\over2}-i\epsilon}r(x)^{{1\over2}+i\epsilon}f(sx)
=q^{{1\over2}-i\epsilon}(s^1f)(x).
\leqno(\a.\f)$$
Suppose $x\in \partial s^+$. If $i>1$, then $F(s_i^{-1}y)=f(x)$ and $((s_i)_*\mu)(y)=\mu(y)$.
If $i=1$, then $F(s_i^{-1}y)=f(sx)$ and $((s_i)_*\mu)(y)=\mu_1(y)=qr(x)\mu(y)$.
Therefore
$$
(s^qf)(x)
=(q-1)f(x)+f(sx)q^{{1\over2}+i\epsilon}r(x)^{{1\over2}+i\epsilon}
=(q-1)f(x)+q^{{1\over2}+i\epsilon}(s^1f)(x).
\leqno(\a.\f)$$
\hfill{$\diamond$}
\smallskip

We would like to extend Lemma \a.1 to a formula expressing the Hecke action as a combination of
Weyl actions for general $w\in W$, in analogy to Proposition 5.1. The appropriate statement
is given in Theorem \a.2. The argument mimics the proof of Proposition 5.1:
we consider the sums indexed by the same posets and split them into sub-sums
in the same way. The differences are: the summands are slightly different;
 $u\in W$ has to be replaced by $x\in\partial W$ (one can think of $x$ as of a limit of $u$).
As before, $P(A|D)$ denotes the set of walls (in $W$) separating $A$ from $D$, but now $A$ and $D$ can be sets of
points in $\partial W$ or elements of $W$. Recall that there is a poset structure on $P(A|D)$. We
will perform summations indexed by anti-chains in this kind of posets. Recall that for $h\in\A P(A|D)$
we denote by $h$ also the product of all reflections corresponding to walls in $h$.
We put $\tau(w,x)=[{dw_*\mu\over d\mu}(x)]^{{1\over2}+i\epsilon}$.

\smallskip
\bf Theorem \a.\t\sl\par
Let $f\in L^2(\partial W,l)$ and let $\epsilon=0$. Then
$$(w^qf)(x)=\sum_ {h\in\A P(1|x,w)}(q-1)^{\# h}q^{(\ell(w)-\# h)/2}((hw)^1f)(x).\leqno(\a.\f)$$

\par
\rm Proof.
We proceed by induction on $\ell(w)$, the first step being given by Lemma \a.1.

Case 1: $1$, $w$ and $x$ are all on the same side of the wall $s$ (the wall separating
$1$ from $s$). Then
$$\eqalign{
s^q(w^qf)(x)
&=q^{1\over2}\tau(s,x)(w^qf)(sx)\cr
&=q^{1\over2}\tau(s,x)\sum_{h \in\A P(1|sx,w)}(q-1)^{\# h}q^{(\ell(w)-\# h)/2}((hw)^1f)(sx)\cr
&=q^{1\over2}\tau(s,x)\sum_{h \in\A P(1|sx,w)}(q-1)^{\# h}q^{(\ell(w)-\# h)/2}\tau(hw,sx)f(w^{-1}h^{-1}sx)\cr
&=q^{1\over2}\sum_{h \in\A P(1|sx,w)}(q-1)^{\# h}q^{(\ell(w)-\# h)/2}\tau(shw,x)f(w^{-1}hsx),
}\leqno(\a.\f)
$$
while the right hand side of the formula for $((sw)^qf)(x)$ that we are trying to prove is
$$
\eqalign{
&\sum_{h \in\A P(1|x,sw)}(q-1)^{\# h}q^{(\ell(sw)-\# h)/2}((hsw)^1f)(x)
=\cr &\sum_{h \in\A P(1|x,sw)}\tau(hsw,x)(q-1)^{\# h}q^{(\ell(sw)-\# h)/2}f(w^{-1}hsx).}
\leqno(\a.\f)
$$

Now, exactly as in the proof of Proposition 5.1, we observe that: (1) $\ell(sw)=\ell(w)+1$;
(2) $P(1|x,sw)=P(1|sx,w)$; (3) if $h\in\A P(1|sx,w)$, then $hs=sh$.

Case 2: The wall $s$ separates $1$ and $w$ from $x$.
$$
\eqalign{
s^q(w^qf)(x)
=&(q-1)(w^qf)(x)+q^{1/2}\tau(s,x)(w^qf)(sx)\cr
=&(q-1)\sum_{h \in\A P(1|x,w)}(q-1)^{\# h}q^{(\ell(w)-\# h)/2}((hw)^1f)(x)\cr
&+q^{1/2}\tau(s,x)\sum_{h \in\A P(1|sx,w)}(q-1)^{\# h}q^{(\ell(w)-\# h)/2}((hw)^1f)(sx)\cr
=&(q-1)\sum_{h \in\A P(1|x,w)}(q-1)^{\# h}q^{(\ell(w)-\# h)/2}\tau(hw,x)f(w^{-1}hx)\cr
&+q^{1/2}\tau(s,x)\sum_{h \in\A P(1|sx,w)}(q-1)^{\# h}q^{(\ell(w)-\# h)/2}\tau(hw,sx)f(w^{-1}hsx)
}\leqno(\a.\f)
$$
while the postulated expression for $((sw)^qf)(x)$ is
$$\eqalign{
&\sum_{k \in\A P(1|x,sw)}(q-1)^{\# k}q^{(\ell(sw)-\# k)/2}((ksw)^1f)(x)=
\cr &\sum_{k \in\A P(1|x,sw)}(q-1)^{\# k}q^{(\ell(w)+1-\#k)/2}\tau(ksw,x)f(w^{-1}skx).
}\leqno(\a.\f)
$$
We consider the same bijection between the summands of (\a.8) and those of (\a.9)
as when we compared (5.6) with (5.7) ($x$ taking the role of $u$). In the case $s\not\in k\in\A P(1|x,sw)$ the corresponding
$h=k^s$; to compare the summands we observe that $\tau(s,x)\tau(hw,sx)=\tau(shw,x)=\tau(ksw,x)$.
\hfill{$\diamond$}
\smallskip

By the same argument, just keeping track of the extra factor with imaginary exponent, one can prove a similar
formula for $\epsilon\ne0$:
$$(w^qf)(x)=\sum_ {h\in\A P(1|x,w)}(q-1)^{\# h}q^{(\ell(w)-\# h)/2+i\epsilon(2\height(h)-\ell(w)-\# h)}((hw)^1f)(x),\leqno(\a.\f)$$
where $\height(h)$ was defined in the previous section (Def.~5.1.c).
\smallskip
We would like to re-formulate Theorem \a.2 in an $x$-free way.
Recall that for a wall $H$ we have defined $H^+$ as the set of all $w\in W$
that are separated from $1$ by $H$. For a collection $h$ of pairwise
intersecting walls (in symbols: $\cap h\ne\emptyset$) we put
$h^+=\cap\{H^+\mid H\in h\}$. We denote
by $\partial h^+$ the set of all $x\in\partial W$ that are separated
from $1$ by each $H\in h$. The characteristic (indicator) function
of a set $U$ will be denoted $\chi_U$. We now re-state Lemma \a.1
and Theorem \a.2.
\smallskip
\bf Corollary \a.\t. \sl\par
Let $f\in L^2(\partial W,l)$, $s\in S$, $w\in W$, and let $\epsilon=0$. Then
$$s^qf=\sqrt{q} (s^1f)+(q-1)\chi_{\partial s^+}\cdot f\leqno(\a.\f)$$
and
$$w^qf=
\sum_{h|\cap h\ne\emptyset, w\in h^+}
(q-1)^{\# h}q^{(\ell(w)-\# h)/2}\chi_{\partial h^+}\cdot (hw)^1f.\leqno(\a.\f)$$
\rm

\bigskip
\centerline{\bf 7. The general setting of right-angled Hecke algebras and its principal series representations}
\def\a{7}
\m=1
\n=1
\medskip

Consider now the case where $q$ is not constant, but is $s$-dependent.
the formulae of sections 5 and 6 still
hold---with
special reading.
First, in the formulae where $e_s$ is involved, $q$ is $q_s$ (these are (5.1), (5.2),
(6.1) and their direct applications).
Second, in powers like $(q-1)^{\# h}$ or $q^{\ell(w)}$, the base $q-1$ or $q$ is a tuple
of numbers $(q_s-1)_{s\in S}$ or $(q_s)_{s\in S}$. The exponent $\# h$, $\ell(w)$ or
$\height(h)$ always counts the walls in a certain set; the walls have types which
are elements of $S$, thus each exponent can be transformed to a multi-index.
For example, $\ell(w)$ is the number of walls separating $1$ and $w$; it can be read
as $(\ell_s(w))_{s\in S}$, where $\ell_s(w)$ is the number of walls of type $s$
that separate $1$ and $w$. Then $q^{\ell(w)}$ is $\prod_{s\in S}q_s^{\ell_s(w)}$;
similarly for other powers.

Next we want to suggest another interpretation of \S6.
We consider the algebra $\H$ defined in \S5, for a given set of complex parameters $q_s$ and
interpretation of the formulae as described above.
One may now read equation (6.1) in Lemma~6.1 (taking $q=q_s$) as a definition of a 
representation of the generators of $\H$ as operators on $L^2(\partial W)$.
It is easy to check by direct calculation that these operators satisfy $(s^q)^2=(q-1)s^q+qI$.
One then reads Theorem~6.2 and its proof as stating that the above gives a well defined representation of $\H$.

\bigskip
\centerline{\bf III. Hecke irreducibility}
\bigskip

\m=1

Our goal in this part is to show that the representation of $\H$ on $L^2(\partial W)$ described in equation~(6.5) is irreducible.
Recall that this will also prove the irreducibility of the representation $\rho$ of $G'$ considered in part I.
Our proof is a modification
of the proof given in [BM]. The main tool which allows us to use [BM]
almost verbatim is the following pointwise inequality between functions on $\partial W$:
$$w^11\le q^{-\ell(w)/2}w^q1\le Cw^11.$$
This inequality is established is section 9. In section 10 we explain how to use it
to prove Hecke irreducibility.

The Hecke algebra representation on $L^2(\partial W)$ is given by equation~(6.5).
As explained in \S7, this formula makes sense for all complex values of $q$,
not just for positive integers. Our proof of irreducibility will work for
all real $q_s\ge1$.
We assume this is the case from now on.
For readability purposes we adopt the conventions explained in \S7, writing our formulae as if $q_s$ was a unique parameter $q$.

The notation for this part is as follows. By $1$ we denote the constant function with value 1 on $\partial W$;
but we also use $1$ for the unit element of the group $W$.
We identify elements of $W$ with points in $\H^d$
via the orbit map $w\mapsto w0$ (hence, $1$ often stands for the point $0$). Distance of $x$ and $y$ in $\H^d$ is denoted $|x-y|$;
in particular,
$|w-w'|$ is the hyperbolic distance from $w0$ to $w'0$. We shorten $|w-1|$ to $|w|$.
The length-parameterized geodesic from $a$ to $b$ will be denoted $\gamma_a^b$;
we shorten $\gamma_0^b$ to $\gamma^b$---thus $\gamma^w=\gamma_1^w=\gamma_0^{w0}$.
The ball of center $m$ and radius $r$ in a metric space $M$ is $M(m,r)$.
The closed $r$-neighborhood of a set $A$ will be denoted $A[r]$.
By $h$ we denote a finite set of pairwise intersecting tessellation walls in $\H^d$, as well
as the composition of reflection across these walls (the group $W$  is right-angled, hence the order of composition
is not important). The intersection of the walls in $h$ is a totally geodesic subspace $\cap h$. The map
$h$ is an isometry of $\H^d$; in fact, it is the orthogonal reflection in $\cap h$.
\bigskip

\centerline{\bf 8. Individual estimates}
\medskip
\def\a{8}
\m=1
\n=1
The goal of this section is to establish estimates
for $((hw)^11)(z)=\tau(hw,z)$. Throughout this section we assume that
$h$ is a collection of pairwise perpendicular tessellation walls,
$w\in h^+$ and $z\in\partial h^+$.
Recall that $\tau(hw,z)=\exp(-\eta\beta_z(1,hw)/2)$, where
$\eta=\d-1$ is the dimension of the boundary $\partial W$ (cf [BM], p.~52). The Busemann function
$\beta_z(1,hw)=\lim_{x\to z}(|1-hw|-2(x|hw))$, the limit being taken over
points $x\in\H^\d$ converging to $z$.
\smallskip
\bf Definition \a.\t\rm\par
For $w\in h^+$ we put
$$(w|\cap h)=|w|-|w-hw|/2.\leqno(\a.\f)$$
\smallskip
This quantity will be investigated (and the notation explained) in the next section.
The estimate for $\tau(hw,z)$ that we need is as follows.
\medskip
\bf Proposition \a.\t\sl\par
Suppose that $h$ is a collection of pairwise perpendicular tessellation walls.
Assume that $w\in h^+$, $z\in\partial h^+$. Then
$$\tau(hw,z)\le e^{\eta\delta}\exp(\eta\min\{(z|w),(w|\cap h)\})\exp(-\eta|w|/2).\leqno(\a.\f)$$
\medskip\rm
This estimate is equivalent to:
$$\beta_z(1,hw)\ge |w|- 2\min\{(z|w),(w|\cap h)\}-2\delta,\leqno(\a.\f)$$
which is the limit (as $x\to z$) of the following pair of inequalities:
\medskip
\bf Proposition  \a.\t\sl\par
Suppose $w,x\in h^+$. Then
\item{a)}{$\displaystyle2(x|hw)-|hw|\le2(x|w)-|w|;$}
\item{b)}{$\displaystyle2(x|hw)-|hw|\le2(w|\cap h)-|w|+2\delta.$}
\medskip\rm
Our main tool will be:
\medskip
\bf Lemma \a.\t\sl\par
Let $h$ be a collection of pairwise perpendicular hyperplanes in $\H^\d$.
Suppose that two points $a,b\in\H^\d$ are not separated by any of these
hyperplanes. Then
$|a-b|\le|a-hb|$.\rm\par
Proof. The geodesic segment $[a,hb]$ intersects every hyperplane of $h$.
Let the intersection points be $p_1,\ldots,p_j$, numbered from $a$ towards $hb$,
and let $h_i$ be the hyperplane passing through $p_i$. Then
$[a,p_1]\cup h_1[p_1,p_2]\cup h_1h_2[p_2,p_3]\cup\ldots\cup h[p_j,hb]$
is a piecewise geodesic path from $a$ to $b$ of length $|a-hb|$.\hfill{$\diamond$}
\smallskip
\bf Proof of Prop.~\a.3.a):\rm \par
We apply Lemma \a.4 to $h$ and the points $x,w\in h^+$ to get $|x-w|\le|x-hw|$. Then
$$2(x|hw)-|hw|=|x|-|x-hw|\le |x|-|x-w|=2(x|w)-|w|.\leqno(\a.\f)$$
\hfill{$\diamond$}
\smallskip
\bf Proof of Prop.~\a.3.b):\rm \par
Expanding the Gromov product and Def.~\a.1 we see that the required inequality is equivalent to
$$|x|-|x-hw|\le|w|-|w-hw|+2\delta,\leqno(\a.\f)$$
which we rewrite as
$$|x-1|+|w-hw|\le |w-1|+|x-hw|+2\delta.\leqno(\a.\f)$$
A basic property of quadrangles in hyperbolic spaces ([BS, sec.~2.4.1]) applied to
$(x,w,1,hw)$ yields
$$|x-1|+|w-hw|\le\max\{|x-w|+|1-hw|,|x-hw|+|w-1|\}+2\delta.\leqno(\a.\f)$$
By Lemma \a.4 we get $|x-w|\le|x-hw|$, $|1-hw|\le|1-w|$. Now (\a.6) follows from (\a.7).
\hfill{$\diamond$}

\bigskip
\centerline{\bf 9. Hecke action estimate}
\def\a{9}
\m=1
\n=1
\medskip
We begin with an explanation of $(w|\cap h)$, and then show the pointwise
inequality $q^{-\ell(w)/2}w^q1\le Cw^11$. Informally,
if $w\in h^+$, then the geodesic from $1$ to $w$ gets close to $\cap h$
at time $|h|/2$ and stays near $\cap h$ till time $(w|\cap h)=|w|-|w-hw|/2$.
\medskip
\bf Lemma \a.\t\sl
$$(w|\cap h)\ge0.$$
\rm Proof. By Lemma 8.4 we have $|1-hw|\le |1-w|$. Using this and the triangle
inequality we get:
$$2(w|\cap h)=2|w|-|w-hw|=|w-1|+|1-w|-|w-hw|\ge
|w-1|+|1-hw|-|w-hw|\ge0.\leqno(\a.\f)$$
\hfill{$\diamond$}

\bf Lemma \a.\t\sl\par
If $x\in h^+$, then $(x|h)\ge|h|/2$.
If $z\in \partial h^+$, then $(z|h)\ge|h|/2$.\rm\par\rm\par
Proof.
By Lemma 8.4, we have $|x|=|x-1|\ge|x-h|$, hence
$2(x|h)=|x|+|h|-|x-h|\ge|h|.$
The boundary version is obtained by passing to limits.\hfill{$\diamond$}
\medskip
\bf Proposition \a.\t \sl\par
Suppose that $w\in h^+$. Then
$\gamma^w([|h|/2,(w|\cap h)])$ is contained in the (closed) $\delta$-neighborhood $\cap h[\delta]$
of $\cap h$.\rm\par
Proof. Since $\cap h[\delta]$ is convex, it suffices to show that
$\gamma^w(|h|/2),\gamma^w((w|\cap h))\in \cap h[\delta]$.

By Lemma \a.2, $(w|h)\ge|h|/2$.
This implies that $|\gamma^w(|h|/2)-\gamma^h(|h|/2)|\le\delta$ (cf.~[BS, def. 1.2.2]).
Now $\gamma^h(|h|/2)$ is the midpoint of $[1,h]$, and it belongs to $\cap h$.

We can now repeat the argument with the points $1$ and $w$ interchanged. The point $hw$ plays the role of the point $h$,
and the distance from $w$ to $\cap h$ is $|w-hw|/2$.
We get $(1|hw)_w\ge |w-hw|/2$, hence $|\gamma_w^1(|w-hw|/2)-\gamma_w^{hw}(|w-hw|/2)|\le\delta$.
Now $\gamma_w^1(|w-hw|/2)=\gamma^w((w|\cap h))$, while $\gamma_w^{hw}(|w-hw|/2)$ is the midpoint
of $[w,hw]$ and belongs to $\cap h$.\hfill{$\diamond$}
\smallskip
\bf Corollary \a.\t\sl\par
Let $w\in h^+$, $z\in \partial h^+$. Then
$\displaystyle\gamma^w(\min\{(z|w),(w|\cap h)\})\in\cap h[2\delta].$

\rm Proof.
If the minimum is $(w|\cap h)$, the claim follows directly from Prop.~\a.3.

Suppose then that $(z|w)<(w|\cap h)$. By hyperbolicity ([BS, def.2.1.6]) and Lemma~\a.2 we have
$$(z|w)\ge \min\{(z|h),(w|h)\}-\delta\ge |h|/2-\delta.\leqno(\a.\f)$$
Therefore
$(z|w)\in[|h|/2-\delta,(w|\cap h)]$, and
$$\gamma^w((z|w))\in
\gamma^w([|h|/2-\delta,(w|\cap h)])\subseteq
\gamma^w([|h|/2,(w|\cap h)])[\delta]\subseteq
\cap h[\delta][\delta]=\cap h[2\delta].\leqno(\a.\f)$$
\hfill{$\diamond$}
\smallskip
To summaries our discussion, we put $s(h)=\min\{(z|w),(w|\cap h)\}$.
(We suppress the dependence on $z$ and $w$ since these will be fixed.)
Then we put together Cor.~\a.4 and Prop.~8.2:
\smallskip
\bf Corollary \a.\t\sl\par
Let $h$ be such that $w\in h^+$, $z\in \partial h^+$. Then
$$\tau(hw,z)\le e^{\eta\delta}\exp(\eta s(h))\exp(-\eta|w|/2).\leqno(\a.\f)$$
Furthermore, $\cap h$ intersects the closed ball of radius $2\delta$
centered at $\gamma^w(s(h))$.\rm
\smallskip
To proceed, we need an estimate of the number of $h$ with a given $s(h)$.
\smallskip
\bf Proposition \a.\t\sl\par
There exists a constant $M$ (depending only on the tessellation)
such that for any $z$ and $w$, and any
(closed) interval $I$ of length $\le 1$ contained in $[0,(z|w)]$,
the number of $h$ satisfying $w\in h^+$, $z\in \partial h^+$, $s(h)\in I$
is at most $M$.\rm
\par
Proof.
For any $h$ as in the statement the set $\cap h$ intersects $(\gamma^w(I))[2\delta]$
(by Corollary \a.5).
The latter set is contained in $\H^d(\gamma^w(\max(I)),1+3\delta)$. Any tessellation chamber intersecting
this ball is contained in $S=\H^d(\gamma^w(\max(I)),1+3\delta+\diam(P))$
(where $P$ is a chamber).
Now all chambers have the same volume, and the volume of a ball is finite and depends only on the radius.
Therefore the number of chambers contained in $S$ is uniformly bounded.
Consequently, the total number of faces of these chambers is also uniformly bounded,
and the latter number is not smaller than the number of $h$'s that we are after. \hfill{$\diamond$}
\smallskip
Note that $s(h)$ always belongs to $[0,(z|w)]$. Indeed, $(w|\cap h)\ge0$ by Lemma \a.1,
while a Gromov product is non-negative by the triangle inequality.
\medskip
\bf Proposition \a.\t\sl\par
There exists a constant $C$ depending only on the tessellation and on
the parameter $q$,
such that for any $w\in W$ we have a pointwise inequality
$$w^11\le q^{-\ell(w)/2}w^q1\le Cw^11$$
\rm\par
Proof.
By Corollary 6.3 we have
$$q^{-\ell(w)/2}w^qf=w^1f+\sum_{j=1}^d\left(q-1\over\sqrt{q}\right)^j
\sum_{h|\# h=j,\cap h\ne\emptyset, w\in h^+}
\chi_{\partial h^+}\cdot (hw)^1f.\leqno(\a.\f)$$
In  particular, for a $z\in \partial W$,
$$(q^{-\ell(w)/2}w^q1)(z)=(w^11)(z)+\sum_{j=1}^d\left(q-1\over\sqrt{q}\right)^j\sum_{h|\# h=j,\cap h\ne\emptyset, w\in h^+,z\in\partial h^+}
\tau(hw,z),\leqno(\a.\f)$$
Let $Q=\max\{({q-1\over\sqrt{q}})^j\mid j=1,\ldots,\d\}$.
In the calculation below $h$ always satisfies $w\in h^+$, $z\in \partial h^+$,
$h\ne\emptyset$---we only state
explicitly the extra conditions. We use Cor.~\a.5.
$$\eqalign{q^{-\ell(w)/2}w^q1(z)&\le w^11(z)+Q\sum_{h}\tau(hw,z)\cr
&\le w^11(z)+Q\sum_he^{2\delta}e^{\eta s(h)}e^{-\eta|w|/2}\cr
&\le w^11(z)+Qe^{2\delta}e^{-\eta|w|/2}\sum_he^{\eta s(h)}
}\leqno(\a.\f)$$
We now focus on the sum:
$$\eqalign{\sum_he^{\eta s(h)}&\le
\sum_{i=0}^{\lfloor(z|w)\rfloor}
\sum_{h\mid s(h)\in[(z|w)-(i+1),(z|w)-i]}e^{\eta s(h)}\cr
&\le
\sum_{i=0}^{\lfloor(z|w)\rfloor}
\#\{h\mid s(h)\in[(z|w)-(i+1),(z|w)-i]\}e^{\eta ((z|w)-i)}\cr
&\le
e^{\eta(z|w)}\sum_{i=0}^{\lfloor(z|w)\rfloor} Me^{-\eta i}
\le
Me^{\eta(z|w)}\sum_{i=0}^\infty e^{-\eta i}\cr
&\le
{M\over 1-e^{-\eta}}e^{\eta(z|w)}.
}\leqno(\a.\f)$$
Recalling that $w^11(z)=\tau(w,z)=e^{\eta((z|w)-|w|/2)}$ and putting
$C=1+e^{2\delta}{QM\over 1-e^{-\eta}}$ we obtain the claim.

In the multi-parameter case, the index $j$ in formulae $(\a.5)$, $(\a.6)$, and in the definition of $Q$,
should be read as a multi-index (with 0-1 components and total degree $\le d$).
  \hfill{$\diamond$}

\bigskip
\centerline{\bf 10. The [BM] argument}
\medskip
\def\a{10}
\m=1
\n=1

Our goal now is to prove irreducibility of the representation $\rho_0$
of the Hecke algebra $\H=H(G,K)$ on the space $V^K$.
The argument follows very closely an irreducibility argument
in [BM]. For easier comparison, we adjust our notation to match [BM].
Thus, the boundary $\partial W$ will be denoted $B$, and the measure
$l$ will be called $\nu$. We also put $H=V^K=L^2(B,\nu)$. The base-point $0$ will be denoted $p$.
The characteristic (indicator) function of a set $U$ is $\chi_U$.
Convergence of operators means weak convergence.

Let $\Gamma$ be a torsion-free finite-index subgroup of $W$.
For $q=1$, the Hecke representation reduces to a group algebra
representation of $W$; this can further be restricted to $\Gamma$.
It is shown in [BM] that that representation of $\Gamma$ is irreducible.
This is achieved by showing that the von Neumann algebra generated
by the representation operators is the whole ${\rm End}(H)$. More precisely,
for a measurable $U\subseteq B$ with zero-measure boundary the operators
$$T^{\chi_U}_t={1\over |S_t|}\sum_{\gamma\in\Gamma}{\chi_U(z(\gamma))\over\den1}\gamma^1\leqno(\a.\f)$$
are shown to converge to $\langle\cdot,1\rangle\chi_U$ as $t\to\infty$.
We explain the notation: $z(\gamma)$ is a short-hand for $z^{\gamma p}_p$---the limit point
of the geodesic ray from $p$ through $\gamma p$. The set $S_t$ is a spherical layer
of $\Gamma$:
$$S_t=\{\gamma\in\Gamma\mid |p-\gamma p|\in(t-R,t+R)\}\leqno(\a.\f)$$
for $R={\rm diam}(\Gamma\setminus \H^\d)$.
In our setting, a completely analogous result is true:
\smallskip
\bf Lemma \a.\t.\sl\par
For any $U\subseteq B$ with zero-measure boundary,
$\langle\cdot, 1\rangle\chi_U$ is a limit point (as $t\to\infty$) of the operators
$${}^qT^{\chi_U}_t=
{1\over |S_t|}\sum_{\gamma\in\Gamma}{\chi_U(z(\gamma))\over\den q}\gamma^q.\leqno(\a.\f)$$

\rm\par
Proof.
A) First, we wish to argue that (for fixed $q$ and $U$, and varying $t$)
the operators ${}^qT^{\chi_U}_t$ have uniformly bounded norms on $L^2(B)$.

These operators map non-negative functions to non-negative functions. Moreover,
for every non-negative $f\in L^2(B)$ we have ${}^qT^{\chi_U}_tf\le{}^qT^{1}_tf$, so that
it is enough to establish a uniform bound for $\|{}^qT^{1}_t\|_{L^2\to L^2}$.
This is done by first estimating
$$\|{}^qT^{1}_t\|_{L^\infty\to L^\infty}=\|{}^qT^{1}_t1\|_{L^\infty}=\left\|
{1\over |S_t|}\sum_{\gamma\in S_t}{\gamma^q1\over\den q}\right\|_{L^\infty}.\leqno(\a.\f)$$
From the pointwise estimate (Prop.~9.7) we get
$q^{-\ell(\gamma)/2}\|\gamma^q1\|_{L^\infty}\le C\|\gamma^11\|_{L^\infty}$,
$q^{-\ell(\gamma)/2}\den q \ge \den1$, hence
$$\|{}^qT^{1}_t\|_{L^\infty\to L^\infty}\le C \left\|
{1\over |S_t|}\sum_{\gamma\in S_t}{\gamma^11\over\den 1}\right\|_{L^\infty}.\leqno(\a.\f)$$
The latter is uniformly bounded by Prop.~4.5 of [BM]
(the notation there is $\lambda^{\gamma p}=\gamma^11$). Next, we observe that
${}^qT^{1}_t$ are self adjoint: for  a generator $s$ of $W$ we have $(s^q)^*=s^q$
(as we see from (6.11)),
hence for $w\in W$ we get $(w^q)^*=(w^{-1})^q$; the set $S_t$ is invariant under taking
inverses, because $|w|=|w^{-1}|$. Consequently, the uniform $L^\infty$--operator norm bound
yields a uniform $L^1$--operator norm bound. Finally, Riesz--Thorin interpolation
gives a uniform $L^2$--operator norm bound.

Banach--Alaoglu theorem implies that the family ${}^qT^{\chi_U}_t$
has a limit point as $t\to\infty$.
Let ${}^qT^{\chi_U}_\infty$ be such a point.

\smallskip

\item{B)} We will show that ${}^qT^{\chi_U}_\infty=\langle\cdot,1\rangle\chi_U$

As in [BM], for $U\subseteq B$ and $a>0$ we denote by $U(a)$ the $e^{-a}$-neighborhood of $U$ in $B$.

\bf Claim \a.\t. \rm (cf Lemma 5.2 in [BM]) \sl\par
Let $V\subseteq B$, $a>0$. There is a constant $C_0$ such that for every $\gamma\in\Gamma$ satisfying
$z(\gamma)\not\in V(a)$ we have
$${\langle\gamma^q1,\chi_V\rangle\over \den q}\le{C_0e^a\over |\gamma|}.\leqno(\a.\f)$$
\rm Proof. For $q=1$, the claim is a part of Lemma 5.2 in [BM]. The general case follows from this special case,
because  $q^{-\ell(\gamma)/2}\gamma^q1\le C\gamma^11$, $q^{-\ell(\gamma)/2}\den q \ge \den1$. {$\diamond$}
\smallskip
\bf Claim \a.\t. \rm (cf Prop.~5.1 in [BM]) \sl\par
Assume we are given a family of elements  $\psi_t\in {\bf R}\Gamma$ with non-negative
coefficients, where $t$ is a real positive parameter. Suppose that for every $t$ we have
$\|\psi_t\|_{L^1}\le 1$, and that for every $\gamma\in\Gamma$ we have
$\lim_{t\to\infty}\psi_t(\gamma)=0$. Then, for every measurable $V\subseteq B$ and every $a>0$,
$$\limsup_{t\to\infty}\sum_{\gamma\in\Gamma}\psi_t(\gamma){\langle\gamma^q1,\chi_V\rangle\over
\den q}\le\limsup_{t\to\infty}\sum_{\gamma\in\Gamma}\psi_t(\gamma)\chi_{V(a)}(z(\gamma))\leqno(\a.\f)$$
\rm\par
Proof. The proof is as in [BM, p.~59], with $\rho(\gamma)$ replaced with $\gamma^q$ everywhere.
An outline is as follows. The sum on the left hand side is divided into three parts:
\item{1)} $|\gamma|<t_0$---this part goes to 0 since $\psi_t\to0$ pointwise;
\item{2)} $|\gamma|>t_0$ and $z(\gamma)\in V(a)$---this is bounded by the right hand side;
\item{3)} $|\gamma|>t_0$ and $z(\gamma)\not\in V(a)$---this is negligible by Claim~\a.2.
{$\diamond$}
\smallskip
Specializing to $\psi_t={1\over |S_t|}\sum_{\gamma\in S_t} \chi_U(z(\gamma))\gamma$ we get
\smallskip
\bf Claim \a.\t. \rm (cf Cor.~5.3 in [BM]) \sl\par
For measurable $U,V\subseteq B$ that are positive distance apart
$$\lim_{t\to\infty}\langle {}^qT^{\chi_U}_t1,\chi_V\rangle=0.\leqno(\a.\f)$$
\smallskip
\rm An easy consequence is
\smallskip
\bf Claim \a.\t. \rm (implicit in [BM]) \sl\par
For measurable $U,V\subseteq B$ that are positive distance apart, and a measurable $W\subseteq B$,
$$\lim_{t\to\infty}\langle {}^qT^{\chi_U}_t\chi_W,\chi_V\rangle=0.\leqno(\a.\f)$$

\rm Proof. $0\le \langle {}^qT^{\chi_U}_t\chi_W,\chi_V\rangle\le
\langle{}^qT^{\chi_U}_t1,\chi_V\rangle\to 0.$ {$\diamond$}
\smallskip
\rm Specializing Claim \a.3 to $\psi_t={1\over |S_t|}\sum_{\gamma\in S_t}
\chi_U(z(\gamma))\gamma^{-1}$ we get
\smallskip
\bf Claim \a.\t. \rm (cf Cor.~5.4 in [BM]) \sl\par
For measurable $U,V\subseteq B$ and every $a>0$
$$\limsup_{t\to\infty}\langle {}^qT^{\chi_U}_t\chi_V,1\rangle\le
\limsup_{t\to\infty}{1\over |S_t|}\sum_{\gamma\in S_t} \chi_U(z(\gamma^{-1}))\chi_{V(a)}(z(\gamma))
\leqno(\a.\f)$$
\smallskip
\rm Using the last two claims we calculate:
$$
\eqalign{
\limsup_{t\to\infty}\langle{}^qT^{\chi_U}_t\chi_V,\chi_W\rangle &=
\limsup_{t\to\infty}\langle{}^qT^{\chi_{U\cap W(a)}}_t\chi_V,\chi_W\rangle+
\limsup_{t\to\infty}\langle{}^qT^{\chi_{U\cap W(a)^c}}_t\chi_V,\chi_W\rangle
\cr
&\le\limsup_{t\to\infty}\langle{}^qT^{\chi_{U\cap W(a)}}_t\chi_V,\chi_{W(2a)}\rangle+0\cr
&\le\limsup_{t\to\infty}\langle{}^qT^{\chi_{U\cap W(a)}}_t\chi_V,1\rangle-
\lim_{t\to\infty}\langle{}^qT^{\chi_{U\cap W(a)}}_t\chi_V,\chi_{W(2a)^c}\rangle\cr
&\le \limsup_{t\to\infty}{1\over |S_t|}\sum_{\gamma\in S_t} \chi_{U\cap W(a)}(z(\gamma^{-1}))
\chi_{V(a')}(z(\gamma))
}
\leqno(\a.\f)
$$
By Margulis' thesis (cf Cor.~C.2 in [BM]) the latter limit exists (if
$\nu(\partial U)=\nu(\partial W(a))=\nu(\partial V(a'))=0$) and equals
$\nu(U\cap W(a))\nu(V(a'))$. We deduce (assuming
$\nu(\partial U)=\nu(\partial W)=\nu(\partial V)=0$)
$$
\limsup_{t\to\infty}\langle{}^qT^{\chi_U}_t\chi_V,\chi_W\rangle\le\nu(U\cap W)\nu(V)\leqno(\a.\f)$$
Replacing $U$ or $V$ or $W$ by its complement we get seven similar inequalities. Adding them
up [BM] deduce (for $q=1$, but the argument is general):
\smallskip
\bf Claim \a.\t. \rm (cf Prop.~5.5 in [BM]) \sl\par
For every measurable $U,V,W\subseteq B$ with zero-measure boundaries
$$\lim_{t\to\infty}\langle{}^qT^{\chi_U}_t\chi_V,\chi_W\rangle=\nu(U\cap W)\nu(V)\leqno(\a.\f)$$
\rm
Proof. Sketched above. For other details see [BM, p.~60]. {$\diamond$}
\smallskip
As explained further in [BM, p.~60], this equality allows to conclude
that ${}^qT^{\chi_U}_\infty=\langle\cdot,1\rangle\chi_U$.
\hfill{$\diamond$}
\medskip

Thus the von Neumann algebra $VN^q$ generated by $\{\gamma^q\mid \gamma\in\Gamma\}$
contains all operators $\langle\cdot,1\rangle\chi_U$ (assuming $\nu(\partial U)=0$).
Hence, it also contains all $\langle\cdot,\chi_V\rangle\chi_U$, as
this is $(\langle\cdot,1\rangle\chi_U)\circ(\langle\cdot,1\rangle\chi_V)^*$
(again, assuming $\nu(\partial U)=\nu(\partial V)=0$). Now Lemma B.3 in [BM]
implies $VN^q={\rm End}(H)$. Hecke irreducibility of $H$ follows by Schur's Lemma.

\bigskip
\centerline{\bf References}
\medskip

\item{[BM]}  U.~Bader and R.~ Muchnik, \it Boundary unitary
representations--irreducibility and rigidity, \rm J.~Mod.~Dyn.~{\bf 5} (2011),
no.~1, 49--69.

\item{[BC]} M.~Bekka and M.~Cowling, \it Some irreducible unitary
representations of $G(K)$ for a simple algebraic group $G$ over an algebraic
number field $K$, \rm Math.~Z.~{\bf 241} (2002), 731--741.

\item{[Bourb]} N.~Bourbaki,
\it Groupes et alg\`ebres de Lie, Chapitres IV--VI, \rm
Masson, 2nd edition, 1981.
%Actualit\'es Scientifiques et Industrielles, No.~1337 Hermann,
%Paris 1968, 288 pp.

\item{[BS]}  S.~Buyalo and V.~Schroeder, \it Elements of asymptotic geometry,
\rm EMS Monographs in Mathematics.
European Mathematical Society (EMS), Z\"urich, 2007. xii+200 pp.

\item{[CT]} I.~Capdeboscq and A.~Thomas, \it Cocompact lattices in
complete Kac-Moody groups with Weyl group right-angled
or a free product of spherical special subgroups, \rm
Math.~Res.~Lett.~{\bf 20} (2013), 339--358.

\item{[CS]} M.~Cowling and T.~Steger, \it The irreducibility of restrictions
of unitary representations to lattices, \rm J.~Reine Angew.~Math., {\bf 420}
(1991), 85--98.

\item{[Dav]} M.~Davis, \it The geometry and topology of Coxeter groups,
\rm London Mathematical Society Monographs Series, 32;
Princeton University Press, Princeton, NJ, 2008. xvi+584 pp.

\item{[DO]} J.~Dymara and D.~Osajda, \it Boundaries of right-angled hyperbolic buildings, \rm Fund.~Math.~{\bf 197} (2007), 123--165.

\item{[FTP]} A.~Fig\`a-Talamanca and M.A.~Picardello, \it Harmonic analysis
on free groups, \rm Lecture Notes in Pure and Applied Mathematics, vol.~{\bf 87},
Marcel Dekker Inc., New York, 1983.

\item{[FTS]} A.~Fig\`a-Talamanca and T.~Steger, \it Harmonic analysis
for anisotropic random walks on homogeneous trees, \rm Mem.~Amer.~Math.~Soc.,
{\bf 110} (1994), xii+68.

\item{[Gar1]} L.~Garncarek, \it Analogs of principal series representations for Thompson's groups $F$ and $T$, \rm Indiana Univ.~Math.~J.,{\bf 61} (2012),
619--626.

\item{[Gar2]} L.~Garncarek, \it Boundary representations of hyperbolic groups,
\rm preprint, arXiv:1404.0903.

\item{[HM]} R.~Howe and A.~Moy, \it Hecke algebra isomorphisms for $GL(n)$ over a $p$-adic field,
\rm J.~Algebra {\bf 131} (1990), no. 2, 388--424.

\item{[Iw]} N.~Iwahori, \it On the structure of a Hecke ring of a Chevalley group over a finite field,
\rm J.~Fac.~Sci.~Univ. Tokyo Sect.~I {\bf 10} (1964), 215--236.

\item{[IM]} N. Iwahori and H. Matsumoto, \it On some Bruhat decomposition and the structure of the Hecke rings of $p$-adic Chevalley groups,
\rm Inst.~Hautes \'Etudes Sci.~Publ.~Math.~{\bf 25} (1965), 5--48.

\item{[RR]} B.~R\'emy and M.~Ronan, \it  Topological  groups of Kac--Moody
type, right-angled twinnings and their lattices, \rm
Comment.~Math.~Helv.~{\bf 81} (2006), 191--219.
\bigskip

Technion, Haifa

{\it E-mail address}: {\tt uri.bader@gmail.com}
\medskip
Uniwersytet Wroc\l awski

{\it E-mail address}: {\tt dymara@math.uni.wroc.pl}
\bye